\documentclass{amsart}

\usepackage{epic, eepic, amsmath, amssymb, amsthm}

\usepackage{my_defs}

\newtheorem{theorem}{Theorem}
\newtheorem{proposition}{Proposition}
\newtheorem{lemma}{Lemma}
\newtheorem{corollary}{Corollary}

{\theoremstyle{definition}
\newtheorem{definition}{Definition} 
\newtheorem{example}{Example} 
\newtheorem{remark}{Remark} 

}

\begin{document}

\title{A Lie Algebra for Fr\"olicher Groups}
\author{Martin Laubinger}
\maketitle

\begin{abstract} 
Fr\"olicher spaces form a cartesian closed category which contains the category of smooth manifolds as a full subcategory. Therefore, mapping groups such as $C^\infty(M,G)$ or $\Diff(M),$ but also projective limits of Lie groups are in a natural way objects of that category, and group operations are morphisms in the category. We call groups with this property Fr\"olicher groups. One can define tangent spaces to Fr\"olicher spaces, and in the present article we prove that, under a certain additional assumption, the tangent space at the identity of a Fr\"olicher group can be equipped with a Lie bracket. We discuss an example which satisfies the additional assumption.
\end{abstract}

%

\section*{Introduction}

Given a category $\mathcal C$ whose objects are sets together with some additional structure, a natural question is whether the morphism sets $\Hom_{\mathcal C}(X,Y)$ can be equipped with that structure, thus making them objects of $\mathcal C.$ For instance in the category of topological spaces, we can equip $C(X,Y)$ with the compact-open topology. If one restricts to compactly generated Hausdorff spaces and defines the topology on $X \times Y$ appropriately, there is a homeomorphism \begin{equation} 
\label{eq:ExpLaw} C(X,C(Y,Z)) \cong C(X \times Y,Z), 
\end{equation}
which simply maps $\phi:X \rightarrow C(Y,Z)$ to the map $\tilde \phi (x,y) = \phi(x)(y).$ Equation (\ref{eq:ExpLaw}) is called the \emph{exponential law}, and categories in which a form of the exponential law holds are called \emph{cartesian closed}.

There have been efforts to define cartesian closed categories of vector spaces and differentiable mappings, leading to the convenient calculus of Kriegl, Michor and others. We refer to the monograph \cite{Kriegl:1997} and the extensive list of references therein.

On the other hand, there are more abstract approaches to `sets with smooth structure'. For instance, there are the closely related \emph{diffeological spaces} of Souriau and \emph{differential spaces} of Chen.  In the present paper we use a notion of smooth spaces due to Fr\"olicher. The smooth structure of a \emph{Fr\"olicher space} is given by curves $c:\R \rightarrow X$ and functions $f:X \rightarrow \R.$ The sets of curves and functions are required to determine one another by requiring that $f$ is a function if and only if $f \circ c$ is in $C^\infty(\R,\R)$ for all curves, and vice versa, $c$ is a curve if and only if $f \circ c$ is in $C^\infty(\R,\R)$ for all functions $f.$ Fr\"olicher, diffeological and other smooth spaces are compared in a recent preprint by Stacey \cite{Stacey:2008}, and also independently in the author's Ph.D. thesis.  

Both diffeological and Fr\"olicher spaces form cartesian closed categories which contain the category of smooth manifolds as a full subcategory. This fact has lead to recent interest in smooth spaces by researchers working in the field of higher categories. We refer to the paper of Baez and Hoffnung \cite{BaHo:2008} and references therein.

Another application is in the field of infinite-dimensional Lie theory, since mapping groups such as $\Diff(M)$ and $C^\infty(M,G),$ for smooth manifolds $M$ and Lie groups $G,$ are in a natural way groups in these categories. The motivating question for the present paper is: How much of the classical Lie theory can be generalized to the category of Fr\"olicher spaces?

 Tangent spaces to diffeological spaces have been defined by Hector \cite{Hector:1995}, and to Fr\"olicher spaces by Fr\"olicher \cite{Frolicher:1980}. The Fr\"olicher tangent spaces are conceptually simpler, since tangent vectors are represented by equivalence classes of smooth curves, similar to the classical case. Our goal in the present paper is to define the structure of a Lie algebra on the tangent space $T_eG$ at identity of a Fr\"olicher group. This will be possible only if $T_eG$ satisfies a certain additional condition, namely that it is in a natural way in bijection with its own tangent space at 0. We verify this condition in a non-trivial example.

Let us mention that Fr\"olicher groups have been investigated by Ntumba and Batubenge \cite{Ntumba:2002}. They show that for Fr\"olicher groups $G,$ the tangent space $T_eG$ can be identified with the space of left invariant vector fields. In the present article, we identify $T_eG$ with derivations of the algebra of smooth functions on $G,$ and the main task is to define an element $[v,w] \in T_eG$ which corresponds to the commutator of the derivations associated to $v$ and $w.$ Diffeological groups have been investigated by Souriau \cite{Souriau:1985}, Donato and Iglesias-Zenmour \cite{Donato:1985},  Hector and Mac{\'{\i}}as-Virg{\'o}s \cite{Hector:2002} and Leslie \cite{Leslie:2003}.


We briefly outline the content of the following sections. 

In the first section, we define Fr\"olicher spaces and then discuss differentiation theory in topological vector spaces. Topological vector spaces yield a family of examples for Fr\"olicher spaces, and the particular example of $C^\infty(\R,\R)$ plays a role in the definition of a Fr\"olicher structure on $\Hom$-sets. We state the most important properties of the category of Fr\"olicher spaces: Existence of limits and colimits,  and cartesian closedness. Finally, we give examples of Fr\"olicher groups which are \emph{not} Lie groups, thereby motivating the notion of Fr\"olicher groups.

The second section is concerned with the tangent functor for Fr\"olicher spaces. We define tangent spaces at a point, and define a Fr\"olicher structure on the disjoint union of all tangent spaces. If $X$ and $Y$ are Fr\"olicher spaces, we show that $T(X \times Y)$ is isomorphic to $TX \times TY,$ and we show that $TX$ agrees with the classical tangent bundle if $X$ is a smooth finite-dimensional manifold.

The third section is the central part of the present article. If $G$ is a Fr\"olicher group, we write $TG$ and $T^2G$ as products, using the canonical group structure of $TG.$ This is entirely analogous to the classical case. If $\g = T_eG,$ we associate to any two vectors $v,w \in \g$ a certain element $ \xi \in T_0 \g.$ If a canonical map $\g \rightarrow T_0\g$ is bijective, we let $[v,w]$ denote the inverse image of $\xi$ under this map. Then we show that $(v,w) \mapsto [v,w]$ is a Lie bracket. To this end, we associate to each $v \in \g$ a derivation $D_v$ of the algebra of smooth functions on $G.$ The main task consists in proving that $D_{[v,w]} = [D_v,D_w],$ where on the right hand side we have the usual commutator of derivations. From this, and linearity and injectivity of $v \mapsto D_v,$ it follows immediately that $[\cdot, \cdot]$ satisfies the axioms of a Lie algebra. 

In section four we discuss an example. The additive group $\R^J$ with its product Fr\"olicher structure is a Fr\"olicher group. Under a certain restriction on the cardinality of $J$ we can prove that $T_0\R^J$ is isomorphic to $\R^J$ as a Fr\"olicher space. We use a result of Kriegl and Michor, based on a theorem of Mazur, which  shows that functions $f:\R^J \rightarrow \R$ of the Fr\"olicher structure factor through a \emph{countable} product $\R^{J_0}.$ Since $\R^{J_0}$ is a Fr\'echet space, we can use results from the first section in order to conclude the proof.

In an appendix, we sketch the argument of Kriegl and Michor used in section four.

This paper contains the main result of our Ph.D. thesis, written at Louisiana State University. We would like to express our thanks to our advisor Jimmie Lawson. Part of the research leading to this article was conducted while the author was guest at Technische Universit\"at Darmstadt and Universit\"at Wien, funded by the Louisiana Board of Regents Grant LEQSF(2005-07)-ENH-TR-21. Thanks to Robert Perlis for making this research visit possible. My gratitude goes to Professors Neeb and Michor for their hospitality and helpful discussions.

\section{Fr\"olicher Spaces}

First let us recall the definition and basic properties of the category of smooth spaces as introduced by Fr\"olicher \cite{Frolicher:1980}. If $A$ and $B$ are sets, let $\Map(A,B)$ denote the collection of set maps from $A$ to $B.$ 

\begin{definition}
a) A Fr\"olicher structure on a set $X$ is a pair $(C,F),$ where $C \subset \Map(\R,X)$ and $F \subset \Map(X,\R),$ such that the following two conditions are satisfied:
\begin{itemize}
\item $C = \{ c:\R \rightarrow X \;|\; (\forall f \in F) \, f \circ c \in C^\infty(\R,\R) \}$
\item $F = \{ f:X \rightarrow \R \;|\; (\forall c \in C) \, f \circ c \in C^\infty(\R,\R) \}.$
\end{itemize}
The elements of $C$ and $F$ are called \emph{curves} and \emph{functions} of the  structure, or \emph{smooth} curves and functions, respectively. The triple $(X,C,F)$ is called a \emph{Fr\"olicher space}. We will usually omit $C$ and $F$ from the notation, if there is no danger of confusion. \\
b) A \emph{morphism} or \emph{smooth map} between Fr\"olicher spaces $(X,C,F)$ and $(Y,C',F')$ is a set map $\phi:X \rightarrow Y$ such that $\phi \circ c \in C'$ for all $c \in C,$ or equivalently $f \circ \phi \in F$ for all $f \in F'.$ \\
c) Fr\"olicher spaces and smooth maps form a category, which we denote $\Fcat.$
\end{definition}

\begin{remark} \label{rem:Saturation}
Let $X$ be a set. Given any set $A \subset \Map(X,\R),$ we let
\[
\mathcal S(A) = \{ c \in \Map(\R,X) \;|\; (\forall f \in A) \, f \circ c \in C^\infty(\R,\R) \}
\]
be the \emph{saturation} of $A.$ Similarly, one can define a saturation $\mathcal S(B) \subset \Map(X,\R)$ for $B \subset \Map(\R,X).$ The condition for $(C,F)$ to be a Fr\"olicher structure then reads
\[
(C,F) = (\mathcal S(F), \mathcal S(C)).
\]
\end{remark}

The category $\Fcat$ is closed under various category theoretical constructions, for instance inductive and projective limits, quotients and subsets. Before we state some properties of the category $\Fcat,$ we study Fr\"olicher structures on topological vector spaces.


\subsection*{Differentiation in Topological Vector Spaces}

We will use a notion of maps of class $C^k$ which is studied in detail, and in a more general context, by Bertram, Gl\"ockner and Neeb in \cite{BeGlNe:04}. 

\begin{definition}
Let $E$ and $F$ be topological vector spaces and $U \subset E$ open. Let $U^{[1]} = \{(x,v,t) \;|\; x \in U, x+tv \in U \} \subset U \times E \times \R.$ Then we say that a map $f:U \rightarrow F$ is of \emph{class $C^1$ on $U$} if there is a continuous map $f^{[1]}:U^{[1]} \rightarrow F$ such that
\[
f(x+tv) -f(x) = tf^{[1]}(x,t,v)
\]
for all $(x,t,v) \in U^{[1]}.$ In other words, $f^{[1]}$ is a continuous extension of the difference quotient $t^{-1}(f(x+tv)-f(x)).$ We say that $f$ is of class $C^2$ if $f^{[1]}$ is of class $C^1,$ and inductively, $f$ is of class $C^{k+1}$ if it is of class $C^k,$ and $f^{[k]}$ is of class $C^1.$ Finally, $f$ is of class $C^\infty$ if it is of class $C^k$ for all $k \in \N.$ We will frequently say `$f$ is $C^k$' and `$f$ is $C^\infty$' rather than $f$ is of class $C^k$ and of class $C^\infty,$ respectively. We will not use the term \emph{smooth} for $C^\infty$-maps, since it is already used for morphisms of Fr\"olicher spaces.
\end{definition}

\begin{remark} 
The authors of \cite{BeGlNe:04} show that if the range is assumed to be a \emph{locally convex} vector space, their definition of $C^k$-map agrees with the more classical notion of differentiability in the sense of Michal and Bastiani (see \cite[Proposition 7.4]{BeGlNe:04}). Note also that if domain and range are Banach spaces, there is the notion of Fr\'echet differentiability, which however leads to the same notion of $C^\infty$-maps as Michal-Bastiani's definition. Since we will only be interested in $C^\infty$-maps, the definition we use generalizes the more classical notions.
\end{remark}

The following theorem is classical. We will combine it with a more recent result to obtain information on Fr\"olicher structures on metrizable topological vector spaces.

\begin{theorem}[Boman's Theorem]  \label{thm:boman}
Let $f \in \Map(\R^n,\R)$ be so that $f \circ c$ is $C^\infty$ whenever $c:\R \rightarrow \R^n$ is $C^\infty.$ Then $f$ is $C^\infty.$
\end{theorem}
\begin{proof}
See \cite{Boman:1967}.
\end{proof}

As an immediate consequence of Boman's theorem we have the following:

\begin{corollary}  \label{cor:findimmfd}
If $M$ is a smooth finite-dimensional manifold, then  
\[(C^\infty(\R,M), C^\infty(M,\R)) \] is a Fr\"olicher structure on $M.$ The category of smooth finite-dimensional manifolds is a full subcategory of $\Fcat.$ 
\end{corollary}

\begin{theorem} \label{thm:genBoman}
Let $E$ and $F$ be topological $\R$-vector spaces. Let $U \subset E$ be open and $f:U \rightarrow F$ a mapping. If $E$ is metrizable, then $f$ is of class $C^k$ if and only if the composition $f \circ c: \R^{k+1} \rightarrow F$ is of class $C^k$ for every $C^\infty$-map $c:\R^{k+1} \rightarrow U.$
\end{theorem}
\begin{proof}
See \cite[Theorem 12.4]{BeGlNe:04}.
\end{proof}

\begin{corollary} \label{cor:Metr}
Let $E$ be a metrizable topological vector space. Then $f:E \rightarrow \R$ is smooth if and only if for all smooth $c:\R \rightarrow E,$ the composition $f \circ c$ is smooth. In particular, there is a Fr\"olicher structure on $E$ having $C^\infty(E,\R)$ as functions.
\end{corollary}
\begin{proof}
Clearly, $f \circ c$ is smooth if both $f$ and $c$ are. Now suppose that $f$ is smooth along smooth curves. Let $d: \R^{k+1} \rightarrow E$ be smooth. Then by assumption, the composition $f \circ d \circ c$ is smooth for all smooth $c.$ Boman's Theorem implies that the map $f \circ d$ is smooth. By Theorem \ref{thm:genBoman}, it follows that $f$ is of class $C^k.$ Since $k$ was arbitrary, $f$ is smooth.
\end{proof}

There is an important class of topological vector spaces $E$ on which there is a Fr\"olicher structure of the form $(C^\infty(\R,E),F).$

\begin{definition}
A locally convex vector space $E$ is \emph{convenient} if the following holds: 
\emph{
Let $c \in \Map(\R,E)$ be a map such that $\phi \circ c$ is $C^\infty$ for all continuous linear functionals $\phi$ on $E.$ Then $c$ is of class $C^\infty.$}
\end{definition}

This is only one of several equivalent characterizations of convenient locally convex spaces, see \cite[Theorem I.2.14]{Kriegl:1997} for more. 
The following theorem is due to Alfred Fr\"olicher.

\begin{theorem} \label{thm:convFrechet}
Fr\'echet spaces are convenient.
\end{theorem}
\begin{proof}
See \cite[Th\'eor\`eme 1]{Frolicher:1981}.
\end{proof}

\begin{corollary} \label{cor:Frechet}
If $E$ is a Fr\'echet space, then the pair $(C^\infty(\R,E), C^\infty(E,\R))$ is a Fr\"olicher structure on $E.$
\end{corollary}
\begin{proof}
Fr\'echet spaces are metrizable and convenient, so this follows from Corollary
\ref{cor:Metr} and Theorem \ref{thm:convFrechet}.
\end{proof}

\begin{example} \label{ex:cinfty}
The vector space $E = C^\infty(\R,\R)$ is a Fr\'echet space when equipped with the topology of compact convergence of all derivatives. By Corollary \ref{cor:Frechet}, the pair $(C^\infty(\R,E),C^\infty(E,\R))$ is a Fr\"olicher structure on $E.$
\end{example}

\subsection*{Properties of the Category $\Fcat$} 

We summarize some properties of the category $\Fcat.$ First we need to introduce the important concepts of initial and final Fr\"olicher structures, which are analogous to initial and final topologies. 

\begin{definition}
Given Fr\"olicher structures $(C,F),(C',F')$ on a set $X,$ we say that $(C,F)$ is \emph{finer} than $(C',F')$ if $C \subset C'$ or equivalently if $F' \subset F.$ Similarly, $(C,F)$ is \emph{coarser} than $(C',F')$ if $C' \subset C.$
\end{definition}

\begin{lemma} \label{lem:coarsefine}
If $X$ is a set and $\tilde F \subset \Map(X,\R),$ then there is a unique coarsest Fr\"olicher structure $(C,F)$ on $X$ with $\tilde F \subset F.$ It is given by $(\mathcal S(\tilde F), \mathcal S^2(\tilde F)),$ where $\mathcal S$ is the saturation from Remark \ref{rem:Saturation}.
Similarly, given $\tilde C \subset \Map(\R,X),$ the pair $(\mathcal S^2(\tilde C), \mathcal S(\tilde C))$ is the finest Fr\"olicher structure on $X$ for which elements of $\tilde C$ are curves. We say that $(C,F)$ is \emph{generated} by $\tilde C$ or $\tilde F,$ respectively.
\end{lemma}
\begin{proof}
It suffices to show that the given pairs are Fr\"olicher structures, since they are then clearly the coarsest and finest containing $\tilde F$ and $\tilde C,$ respectively.

Let us only prove that $(\mathcal S(\tilde F), \mathcal S^2(\tilde F))$ is a Fr\"olicher structure, the proof for the other pair being analogous.
By Remark \ref{rem:Saturation}, we need to show that
\[
(\mathcal S(\tilde F), \mathcal S^2(\tilde F))=(\mathcal S^3(\tilde F), \mathcal S^2(\tilde F)).
\]
In general we have $\mathcal S(A) \subset \mathcal S^3(A),$ and we want to prove equality. But since we also have $A \subset \mathcal S^2(A),$ this follows easily:
\[
f \in \mathcal S^3(A) \implies \forall a \in A: f \circ a \in C^\infty(\R,\R)
\implies f \in \mathcal S(A).
\]
\end{proof}

Existence of finest and coarsest Fr\"olicher structures containing a given set of maps allows us to make the following definitions.

\begin{definition}
a) Let $X$ be a set, $\{(X_i,C_i,F_i)\}_{i \in I}$ a collection of Fr\"olicher spaces, and $g_i:X_i \rightarrow X$ and $f_i:X \rightarrow X_i$ set maps. The \emph{initial Fr\"olicher structure with respect to the maps $f_i$} is the Fr\"olicher structure generated by all $f \circ f_i$ with $i \in I$ and $f \in F_i.$ Similarly, the \emph{final Fr\"olicher structure with respect to the maps $g_i$} is the Fr\"olicher structure generated by all $g_i \circ c$ with $i \in I$ and $c \in C_i.$

In particular, if $X$ is a Fr\"olicher space and $\iota: A \rightarrow X$ the inclusion of a subset and $\pi:X \rightarrow B$ the projection onto a quotient, then the \emph{subset structure on $A$} is the initial structure with respect to $\iota,$ and the \emph{quotient structure on $B$} is the final structure with respect to $\pi.$

b) We refer to \cite{MacLane:1971} for the precise definition of diagram, limit and colimit in a category. Given a diagram with objects $(X_i,C_i,F_i)$ in the category $\Fcat,$ we can form the limit $X$ and colimit $Y$ of the underlying sets. The limit and colimit come equipped with maps $a_i:X \rightarrow X_i$ and $\alpha_i:X_i \rightarrow Y.$ If we equip the sets $X$ and $Y$ with the initial and final Fr\"olicher structure with respect to the maps $a_i$ and $\alpha_i,$ it is easy to show that one gets limit and colimit in the category of Fr\"olicher spaces. 
\end{definition}

\begin{example}[Product Structure] \label{ex:product}
Consider the discrete diagram with objects $\{ (X_i,C_i,F_i) \}_{i \in I}.$ The limit of the underlying sets is the cartesian product $\prod_{i \in I} X_i,$ and the \emph{product Fr\"olicher structure} is the coarsest for which all projections $\pi_i$ are morphisms. By Lemma \ref{lem:coarsefine}, its curves are precisely the maps $c = (c_i)_{i \in I}: \R \rightarrow \prod_{i \in I} X_i$ such that $c_i \in C_i$ for all $i \in I.$
\end{example}
%
%
\begin{definition} \label{def:MappingStrct}
Let $C^\infty(\R,\R)$ carry the Fr\"olicher structure as in Example \ref{ex:cinfty}. Given Fr\"olicher spaces $(X,C,F)$ and $(Y,C',F'),$ equip the set $\Hom(X,Y)$ with the initial Fr\"olicher structure with respect to the collection $\{ \Phi_{f,c} \;|\; f \in F', c \in C \}$ of maps 
\[
\Phi_{f,c}:\Hom(X,Y) \rightarrow C^\infty(\R,\R), \quad \phi \mapsto f \circ \phi \circ c.
\]
\end{definition}

The following theorem says that the category $\Fcat$ is cartesian closed, and that composition yields smooth maps between $\Hom$-objects.

\begin{theorem} \label{thm:CC}
a) If $\phi \in \Hom(X \times Y,Z),$ we let $\bar \phi \in \Map(X,\Map(Y,Z))$ be defined by $\bar \phi(x)(y) = \phi(x,y).$ If $\phi$ is smooth, then so are all $\bar \phi(x).$ Hence $\bar \phi$ takes values in $\Hom(Y,Z),$ and it is also smooth. The map
\[
\Hom(X \times Y, Z) \rightarrow \Hom(X,\Hom(Y,Z)), \quad \phi \mapsto \bar \phi
\]
is an isomorphism of Fr\"olicher spaces.

b) If $\phi:Y \rightarrow Z$ is smooth, then so is the map
\[
L_\phi: \Hom(X,Y) \rightarrow \Hom(X,Z), \quad f \mapsto \phi \circ f
\]
of left composition with $\phi.$
\end{theorem}
\begin{proof}
Part a) is proven in \cite[Theorem 23.2]{Kriegl:1997}, and part b) is an easy consequence, see \cite[Corollary 3.13]{Kriegl:1997}.
\end{proof}

Together with the fact that manifolds form a full subcategory of $\Fcat, $ we get natural Fr\"olicher structures on the sets $C^\infty(M,N)$ and $\Diff(M)$ for manifolds $M$ and $N.$ A map $\phi$ from a Fr\"olicher space into one of those sets is smooth if and only if $\tilde \phi(x,y) = \phi(x)(y)$ is smooth on $X \times M.$


\subsection*{Examples of Fr\"olicher Groups} \label{subsct:Examples}

Let $G$ be a group with multiplication $m$ and inversion $i,$ and let $(C,F)$ be a Fr\"olicher structure on $G.$ We say that $G$ is a \emph{Fr\"olicher group} if $i$ is smooth, and if $m:G \times G \rightarrow G$ is smooth with respect to the product structure on $G \times G.$

We discuss some classes of examples.

\begin{example}[Mapping Groups]
Let $X$ be a Fr\"olicher space and $G$ a Fr\"olicher group with multiplication $m_G$ and inversion $i_G.$ Equip $\mathcal G= \Hom(X,G)$ with its Fr\"olicher structure from Definition \ref{def:MappingStrct}, and the group structure given by pointwise multiplication and inversion. We have
\[
fg = m_G \circ (f,g)
\]
and 
\[
f^{-1} = i_G \circ f.
\]
Thus, inversion and multiplication in $\mathcal G$ are smooth by Theorem \ref{thm:CC}.
\end{example}

\begin{example}[Diffeomorphism Groups]
Let $M$ be a Banach manifold. Equip $\Diff(M) \subset \Hom(M,M)$ with the subset Fr\"olicher structure. The group multiplication is given by composition and therefore smooth by Theorem \ref{thm:CC}. Smoothness of inversion is implied by the following theorem, which can be found in the forthcoming book by Gl\"ockner and Neeb \cite{GlNe:2010}. For convenience, we will give the proof.
\begin{theorem}
Let $N$ be a locally convex manifold and $M$ a Banach manifold. Let $f: N \rightarrow \Diff(M)$ be a map such that $f_1 = \tilde f:N \times M \rightarrow M$ is $C^\infty$. Then the map $f_2:(p,m) \mapsto f(p)^{-1}(m)$ is also $C^\infty.$
\end{theorem}
\begin{proof}
Assume without loss of generality that $N$ is open in some locally convex space $E.$ Assume that $M$ is modeled on the Banach space $F.$ Fix $(p_0,m_0) \in N \times M$ and let $m_0'=f(p_0)(m_0).$ We can find charts $(\phi,V)$ and $(\phi,V')$ around $m_0$ and $m_0',$ and an open neighborhood $U$ of $p_0$ such that $f_1(U \times V) \subset V'.$  We let
\[
\Phi:= \phi' \circ f_1 \circ (\id_U \times \phi^{-1}),
\]
which is smooth from an open subset of $E \times F$ to $F,$ and if we let $\Phi_{p_0}(x) = \Phi(p_0,x),$ then this is a local diffeomorphism, since $\Phi_{p_0} = \phi' \circ f(p_0) \circ \phi^{-1}.$ In particular, the differential $d_x \Phi_{p_0}$ is invertible for $x \in \phi(V).$ Now Gl\"ockners generalized implicit function theorem \cite[Theorem 2.3]{gloeckner:2006} implies the existence of open neighborhoods $U_1$ of $p_0$ and $V_1$ of $m_0$ such that
\begin{itemize}
\item $W = \cup_{p \in U_1} \{p\} \times \Phi_{p}(\phi(V_1))$ is open in $E \times F.$
\item The map $(p,m) \mapsto (p,\Phi_p(m))$ is a diffeomorphism of $U_1 \times \phi(V_1)$ onto $W.$
\item The map $(p,v) \mapsto \Phi_p^{-1}(v)$ is smooth on $W.$
\end{itemize}
The last item implies that $f_2(p,m) = (\phi^{-1} \circ \Phi_p^{-1} \circ \phi') (p,m)$ is smooth in a neighborhood of $(p_0,m_0').$  
\end{proof}
\end{example}

\begin{example}[Limits of Lie Groups] \label{ex:limits}
It follows from Corollary  \ref{cor:findimmfd} that all finite dimensional Lie groups are Fr\"olicher groups. Given a diagram of Fr\"olicher groups, we can form the limit in the category of groups and in the category of Fr\"olicher spaces. Both limits have the same underlying set $G,$ and we claim that this limit is a Fr\"olicher group. But the Fr\"olicher structure is generated by the projections $\pi_i: G \rightarrow G_i,$ which are also group homomorphisms. This easily implies that the group operations of $G$ are smooth.

In particular, any product $\prod_{ i \in I} G_i$ of Lie groups is a Fr\"olicher group. We will consider the special case $(\R^J,+)$ in Section \ref{sct:Example}.
\end{example}


\section{The Tangent Functor}

The definition of tangent and cotangent space to a Fr\"olicher space goes back to Fr\"olicher \cite{Frolicher:1986}. If $X$ is a Fr\"olicher space, we equip the disjoint union $TX$ of all tangent spaces with a natural Fr\"olicher structure.  We prove functoriality of $X \mapsto TX,$ show that $T(X \times Y)$ is isomorphic to $TX \times TY$ and finally show that $TM$ is diffeomorphic to the `classical' tangent bundle for smooth finite-dimensional manifolds $M.$

\begin{definition}
a) Let $(X,C,F)$ be a Fr\"olicher space and $x \in X.$ We denote by $C_x$ the set of curves satisfying $c(0)=x.$ We introduce equivalence relations on $C_x$ and $F$ as follows:
\begin{itemize}
\item $c \sim d$ iff $(f \circ c)'(0) = (f \circ d)'(0)$ for all $f \in F,$
\item $f \sim_x g$ iff $(f \circ c)'(0) = (g \circ c)'(0)$ for all $c \in C_x.$
\end{itemize}   
b) We define the tangent space and cotangent space to $X$ at $x$ as $T_xX = C_x /\! \sim$ and $T^xX = F/\! \sim_x,$ respectively. Note that $T^xX$ is a quotient vector space of $F.$ We denote the equivalence classes of curves and functions by $[c]$ and $[f]$ respectively. We will frequently use the notation $v = [s \mapsto c(s)]$ for tangent vectors. \\
c) Let $b:T_xX \times T^xX \rightarrow \R$ be given by $b([c],[f]) = (f \circ c)'(0)$. This is well defined by definition of tangent and cotangent space. Clearly, $b$ is linear in the second argument.\\
d) We say that $u = v + sw$ for $v,w \in T_xX$ and $s \in \R$ if $b(u,\xi) = b(v,\xi) + sb(w,\xi)$ holds for all $\xi \in T^xX.$ The chain rule implies that if $v = [c],$ then the curve $t \mapsto c(st)$ represents $sv.$ Thus, we can always multiply elements  of the tangent space by scalars. We say that \emph{$T_xX$ is a vector space} if $v+w$ exists for all tangent vectors $v,w \in T_xX.$
\end{definition}

The following example illustrates that $T_xX$ is not always a vector space.

\begin{example}
If $X \subset \R^2$ is the coordinate cross with its subset Fr\"olicher structure, consider tangent vectors $v, w \in T_{(0,0)}X$ represented by $t \mapsto (t,0)$ and $t \mapsto (0,t)$ respectively. Then there is no vector $u$ with $u = v+w,$ since a curve $c:\R \rightarrow \R^2$ representing $u$ would satisfy $c'(0) = (1,1).$ But such a curve can not map into $X.$
\end{example}

\begin{lemma} \label{lem:vectoraddition}
If $G$ is a Fr\"olicher group, then all tangent spaces $T_gG$ are vector spaces.
In particular, sum and additive inverse in $T_gG$ are given by
\begin{equation} \label{eq:addition} \begin{split}
[c] + [d] &= [s \mapsto c(s)g^{-1}d(s)]   \\
-[c]  &= [ s \mapsto g c(s)^{-1} g]
\end{split} \end{equation}
\end{lemma} 
\begin{proof}
If $[c], [d] \in T_gG,$ consider the smooth function $\gamma:(s,t) \mapsto c(s)g^{-1}d(t)$ from $\R^2$ into $G.$ We get a curve through $g$ after composing with the diagonal $\Delta: \R \rightarrow \R^2.$ Let $f$ be a function on $G.$ The chain rule yields
\[
(f \circ \gamma \circ \Delta)'(0) =  \frac{\partial}{\partial s} (f \circ \gamma)(s,0) + \frac{\partial}{\partial t} (f \circ \gamma)(0,t) 
\]
and, since $\gamma(s,0)=c(s)$ and $\gamma(0,t) = d(t),$ the latter expression equals $(f \circ c)'(0)+ (f \circ d)'(0),$ which proves that $\gamma \circ \Delta$ represents $[c]+[d].$ 
To check the formula for $-[c],$ just note that applying the sum formula one sees that $[c]+(-[c])$ is represented by
\[
s \mapsto c(s)g^{-1}gc(s)^{-1}g = g.
\]
This is the constant curve through $g,$ hence $[c]+(-[c])$ is the zero vector in $T_gG.$
\end{proof}

\begin{definition}
The \emph{tangent bundle $TX$} of a Fr\"olicher space $X$ is the 
disjoint union of all $T_xX$ for $x \in X.$ If $\phi: X \rightarrow Y$ is a smooth map between Fr\"olicher spaces, we define $T\phi: TX \rightarrow TY$ by $T\phi([c]) = [\phi \circ c].$ 
\end{definition}

\begin{lemma}
Let $\phi:X \rightarrow Y$ be a smooth map of Fr\"olicher spaces. If $T_xX$ is a vector space, then vectors in the image of $T_x\phi$ can be added, and $T_x\phi$ is a linear map onto its image.
\end{lemma}
\begin{proof}
Let $[h] = [c]+[d] \in T_xX.$ Then $T\phi([h]) = [\phi \circ h].$ Let $f$ be smooth on $Y.$ Then $f \circ \phi$ is smooth on $X,$ hence $(f \circ \phi \circ h)'(0) =  (f \circ \phi \circ c)'(0) + (f \circ \phi \circ d)'(0)$ by assumption. This proves additivity, and it is also straightforward to see that $T_x\phi$ commutes with scalar multiplication.
\end{proof}

It is clear that for manifolds, the Fr\"olicher tangent bundle agrees as a set with the classical tangent bundle. In particular, $T\R$ is identified with $\R^2$ via $[c] \mapsto (c(0),c'(0)).$

Given a Fr\"olicher space $(X,C,F),$ we will from now on equip $TX$ with the initial Fr\"olicher structure with respect to all $Tf:TX \rightarrow T\R$ with $f \in F.$

\begin{remark}
One could also use the final structure on $TX$ with respect to the differentials $Tc$ of curves $c \in C.$ Since all $Tf \circ Tc = T(f \circ c)$ are smooth, this would yield a \emph{finer} structure than the initial structure with respect to the $Tf.$
\end{remark}

We now characterize smooth curves into $TX.$

\begin{definition}
Let $c: \R \rightarrow TX$ be a smooth curve. For each $s \in \R,$ pick a representative $c_s: \R \rightarrow X$ of $c(s).$ Then we say that
\begin{equation*} \begin{split} 
\gamma: \R^2 &\rightarrow X    \\
\gamma(s,t) &:= c_s(t)
\end{split} \end{equation*}
\emph{represents} $c.$
The smooth curve $c:\R \rightarrow TX$ yields an element $\xi = [c] \in T^2X.$ We also say that \emph{$\xi$ is represented by $\gamma.$}
\end{definition}

\begin{remark} \label{rem:secondDerivative}
Suppose that $\xi \in T^2X$ is represented by $\gamma: \R^2 \rightarrow X,$ and that $\phi:X \rightarrow Y$ is smooth. Then
\[
T^2\phi(\xi) = [s \mapsto T\phi ([c_s])] = [s \mapsto [t \mapsto \phi(\gamma(s,t))]] 
\]
which shows that $T^2\phi(\xi)$ is simply represented by $\phi \circ \gamma.$
\end{remark}

\begin{lemma} \label{lem:curvesinTX}
A map $\gamma:\R^2 \rightarrow X$ represents a curve $c:\R \rightarrow TX$ if and only if it satisfies the following three conditions:
\begin{itemize}
\item[i)] For fixed $s,$ the map $t \mapsto \gamma(s,t)$ is smooth.
\item[ii)] The map $s \mapsto \gamma(s,0)$ is smooth.
\item[iii)]  The derivative $\frac{\partial (f \circ \gamma)}{\partial t}(s,0)$ is a smooth function in $s$ for all functions $f$ on $X.$
\end{itemize}
\end{lemma}
\begin{proof}
Clearly $i)$ is necessary in order to define $c:\R \rightarrow TX$ by $c(s) = [t \mapsto \gamma(s,t)].$ Now let $f:X \rightarrow \R,$ and consider $Tf \circ c.$ We have
\[
Tf(c(s)) = \left( f(\gamma(s,0)) ,  \frac{\partial (f \circ \gamma)}{\partial t}(s,0) \right)
\]
after identifying $T\R \cong \R^2.$ This shows that $c$ is smooth if and only if $ii)$ and $iii)$ are satisfied.
\end{proof}

\begin{lemma} \label{lem:Tfunc}
Let $X$ and $Y$ be Fr\"olicher spaces, and $\phi:X \rightarrow Y$ a smooth map. Then the projection $\pi:TX \rightarrow X$ and the differential $T\phi:TX \rightarrow TY$ are smooth. In particular, $X \mapsto TX$ defines a functor $T: \Fcat \rightarrow \Fcat.$ 
\end{lemma}
\begin{proof}
Let $\pi:TX \rightarrow X$ and $\pi':T\R \cong \R^2 \rightarrow \R$ be the projections of the tangent bundles. If $f$ is a function on $X,$ then $f \circ \pi = \pi' \circ Tf.$ Hence $f \circ \pi$ is smooth for all such $f,$ therefore $\pi$ is smooth. Now if $\phi:X \rightarrow Y$ is smooth and $f$ is a function on $Y,$ then $T\phi \circ Tf = T(\phi \circ f)$ implies smoothness of $T\phi.$
\end{proof}

\begin{lemma}
For Fr\"olicher spaces $X$ and $Y,$ the spaces $T(X \times Y)$ and $TX \times TY$ are isomorphic.
\end{lemma}
\begin{proof}
Let $\pi_1,\pi_2$ be the projections of $X \times Y$ onto $X$ and $Y,$ respectively, and write $c = (c_1,c_2)$ for curves into $X \times Y.$
The map $\Pi = (T\pi_1,T\pi_2)$ is smooth by the previous lemma. It is onto, because $([c],[d]) \in TX \times TY$ is the image of $[s \mapsto (c(s),d(s))] \in T(X \times Y).$ Now let us show that $\Pi$ is injective. Suppose that $\Pi([c]) = \Pi([d]).$ Then $[c_i]=[d_i]$ for $i = 1,2.$ Let $f$ be a smooth function on $X \times Y.$ Write $c = c_1 \times c_2 \circ \Delta,$ where $\Delta:\R \rightarrow \R^2$ is the diagonal map. Then
\[
(f \circ c)'(0) = (f \circ (c_1 \times c_2) \circ \Delta)'(0) = \frac{\partial (f \circ (c_1,c_2))}{\partial s}(0,0) + \frac{\partial (f \circ (c_1,c_2))}{\partial t}(0,0)
\]
by the chain rule. The derivatives on the right hand side can be rewritten as $(f_1 \circ c_2)'(0)$ and $(f_2 \circ c_1)'(0),$ where $f_1(y)=f(c_1(0),y)$ and $f_2(x)=f(x,c_2(0)).$ The maps $f_1$ and $f_2$ are smooth on $Y$ and $X,$ respectively. Therefore, $(f_1 \circ c_2)'(0) = (f_1 \circ d_2)'(0)$ and similarly for $f_2.$ We conclude
\[
(f \circ c)'(0) = (f_1 \circ d_2)'(0) + (f_2 \circ d_1)'(0) = (f \circ d)'(0),
\]
hence $[c]=[d].$ 
It remains to prove that $\Pi^{-1}$ is smooth. Let $c= (c_1,c_2):\R \rightarrow TX \times TY$ be a curve. If $c_i$ is represented by $\gamma_i,$ then clearly $\Pi^{-1} \circ c$ is represented by $(s,t) \mapsto (\gamma_1(s,t),\gamma_2(s,t)).$ 
We use Lemma \ref{lem:curvesinTX}. Conditions i) and ii) are immediately verified. Property iii) can be proven by a similar trick as injectivity of $\Pi,$ using the chain rule and the diagonal map, as well as the fact that property iii) holds for the $\gamma_i.$
\end{proof}

\subsection{Tangent Functor for Finite Dimensional Manifolds}

The Fr\"olicher tangent functor extends the classical tangent functor. In order to prove this, let us temporarily denote the Fr\"olicher tangent functor by $\tilde T$ and the classical one by $T.$ Furthermore, fix a smooth $n$-dimensional manifold $M$ with charts $(U_i, \phi_i)$ indexed by a set $I.$ Then $\tilde TM$ can be described as set of equivalence classes $[x,i,a]$ where $x \in M, i \in I$ and $a \in \R^n.$ Triples $(x,i,a)$ and $(y,j,b)$ are equivalent if $x=y$ and $d_{\phi_i(x)}(\phi_j \circ \phi_i^{-1})(a)=b.$ The manifold structure of $TM$ is given by the charts $T\phi_i: TU_i \rightarrow U_i \times \R^n, [x,i,a] \mapsto (\phi_i(x),a).$

We define a map $\Phi:\tilde TM \rightarrow TM$ as follows. If $[c] \in \tilde TM,$ let $x = c(0)$ and choose $i \in I$ such that $x \in U_i.$ Now let $a = (\phi_i \circ c)'(0)$ and set $\Phi([c]) = [x,i,a].$

\begin{theorem}
The map $\Phi: \tilde TM \rightarrow TM$ is well defined and a diffeomorphism of Fr\"olicher spaces.
\end{theorem} 
\begin{proof}
Suppose that $[c] = [d].$ Then $c(0)= d(0),$ and we can pick $i \in I$ such that $\Phi([c])=[x,i,a]$ and $\Phi([d]) = [x,i,b].$ It remains to show that $a=b.$ We can use a smooth bump function $\rho$ which is supported in $U_i$ and is constant 1 on a smaller neighborhood of $x,$ to define $f = \rho \phi_i.$ Since $[c] = [d],$ we get $a=(f \circ c)'(0)= (f \circ d)'(0) = b.$ This proves well-definedness of $\Phi.$  

Next we check smoothness of $\Phi.$ Let $c:\R \rightarrow \tilde TM$ be a curve. Then $\pi \circ c:\R \rightarrow M$ is smooth, and we can choose a neighborhood $U$ of $0 \in \R$ which is mapped into a chart $U_i.$ Compose $\Phi \circ c:U \rightarrow TM$ with $T\phi_i$ to get the map
\[
s \mapsto ((\phi_i \circ \pi \circ c)(s), (\phi_i \circ c_s)'(0)),
\]
where $c_s$ represents the vector $c(s) \in \tilde TM.$ This map is smooth by Lemma \ref{lem:curvesinTX}.

Let us construct an inverse map to $\Phi.$ If $U_i \times \R^n$ is a chart for $TM,$ we map $(x,a) \in U_i \times \R^n$ to the vector in $\tilde TM$ represented by the curve $[ t \mapsto \phi_i^{-1}(x + ta)].$ This map is easily seen to be inverse to $\Phi.$ To see that it is smooth, recall that $\tilde TM$ carries the initial Fr\"olicher structure with respect to the maps $Tf$ for all smooth functions $f$ on $M.$ If $v \in T_xX,$ we can represent $v$ by a curve $t \mapsto \gamma(t)$ with $\gamma(0)=x.$ By construction, $\Phi^{-1}(v) = [c],$ and $(Tf \circ \Phi^{-1})(v) = d_xf(v)$ is the classical differential. But $df$ is smooth, which shows that $\Phi^{-1}$ is smooth. 
\end{proof}


\section{The Lie Algebra}

Throughout this section, let $G$ be a Fr\"olicher group and $\g = T_eG.$ 
Let $v = [c]$ and $w=[d]$ be elements of $\g.$ Our strategy is to first define an element $\xi \in T_0\g$ which is related to the commutator of $v$ and $w$ in $T^2G.$ Then we consider a natural map $\g \rightarrow T_0\g$; if this map is bijective, we let $[v,w]$ be the inverse image of $\xi$ under this map. Lastly, we will view elements of $\g$ as derivations on the algebra $F$ of functions on $G$ to prove that $(v,w) \mapsto [v,w]$ is a Lie bracket.

\subsection*{Trivialization of Tangent Bundles}

If $G$ is a Fr\"olicher group with multiplication $m$ and inversion $i,$ then the differentials $Tm$ and $Ti$ are smooth maps, and the chain rule $T(f \circ g) = Tf \circ Tg$ easily implies that $TG$ is a Fr\"olicher group with multiplication $Tm$ and inversion $Ti.$ We simply write $vw$ for $Tm(v,w)$ and note that $[c][d]$ is represented by the curve $s \mapsto c(s)d(s).$ Furthermore, the zero section $G \rightarrow TG$ is a homomorphism, and we write $g[c] = [s \mapsto gc(s)]$ for the left action of $G$ on $TG.$
For classical Lie groups $G$, there are two closely related ways to trivialize $T^2G.$
The first one is to trivialize $TG \cong G \times \g,$ and then apply $T$ to get $T(G \times \g) \cong TG \times \g \times T_0\g,$ then trivialize the first factor.
The second way is to regard $TG$ as a Lie group, trivialize $TTG \cong TG \times T_0(TG)$ and then trivialize $TG.$ See \cite{Didry:2006} for the precise relation of these trivializations.

We will use the first approach. The projection $\pi:TG \rightarrow G$ which maps $[c]$ to $c(0)$ is a group homomorphism. For any $v \in TG$ one has $\pi(v)^{-1}v \in \g.$ 
Define
\[
\Phi: TG \rightarrow G \times \g
\]
by $\Phi([c]) = (c(0),c(0)^{-1}[c]).$

\begin{lemma}
The map $\Phi$ is an isomorphism of Fr\"olicher groups, when $G \times \g$ carries the semidirect product group structure with multiplication
\[
(g,v)(h,w) = (gh, \Ad(h^{-1})v + w).
\]
\end{lemma}
\begin{proof}
Clearly, $\Phi$ is smooth and bijective. The inverse of $\Phi$ is given by $(g,v) \mapsto gv,$ which is also smooth. It remains to show that $\Phi$ is a homomorphism. By definition we have
\[
\Phi([c][d]) = (c(0)d(0), d(0)^{-1}c(0)^{-1}[c][d]).
\]
Let us write $g = c(0)$ and $h = d(0).$ The second component of above expression can be written
\[
h^{-1}g^{-1}[c]hh^{-1}[d].
\]
By \ref{lem:vectoraddition}, this is the sum $h^{-1}g^{-1}[c]h + h^{-1}[d] = \Ad(h^{-1})g^{-1}[c] + h^{-1}[d],$ and $\Phi $ is a homomorphism.
\end{proof}

The map $T\Phi$ is an isomorphism between $T^2G$ and $T(G \times \g).$ 
Further, $T(G \times \g)$ is isomorphic to $TG \times T\g,$ and since $G$ and $\g$ are Fr\"olicher groups, we can further trivialize to get an isomorphism, which we also denote $\Phi,$ between $T^2G$ and $G \times \g \times \g \times T_0\g.$ If $\xi \in T^2G,$ let 
\[
\Phi(\xi)= (\pi_1(\xi),\pi_2(\xi),\pi_3(\xi),\pi_4(\xi)).
\]
We now describe the components $\pi_i(\xi).$ 

\begin{lemma}
Suppose $\gamma:\R^2 \rightarrow G$ represents $\xi \in T^2G.$ Let $g = \gamma(0,0).$ Then we have 
\begin{equation*} \begin{split}
\pi_1(\xi) &= g  \\
\pi_2(\xi) &= [s \mapsto g^{-1} \gamma(s,0)] \\
\pi_3(\xi) &= [t \mapsto g^{-1} \gamma(0,t)]  \\
\pi_4(\xi) &= [s \mapsto [t \mapsto \gamma(0,t)^{-1} g \gamma(s,0)^{-1} \gamma(s,t) ]]   \in T_0\g.
\end{split} \end{equation*}
\end{lemma}
\begin{proof}
We first apply $T\Phi:$
\begin{equation*} \begin{split}
T\Phi(v) &= [s \mapsto \Phi([t \mapsto \gamma(s,t)])]   \\
&= [s \mapsto (\gamma(s,0), \gamma(s,0)^{-1}[t \mapsto \gamma(s,t)]) \\
&= [s \mapsto (\gamma(s,0), [t \mapsto \gamma(s,0)^{-1}\gamma(s,t)])].
\end{split} \end{equation*}
We use the canonical identification $T(G \times \g) \cong TG \times T\g.$ Then we use the trivializations of $TG$ and $T\g,$ respectively. Clearly, the part in $TG = G \times \g$ is given by
\[
(\gamma(0,0), [s \mapsto \gamma(0,0)^{-1}\gamma(s,0)]),
\]
which proves the formula for $\pi_1$ and $\pi_2.$ Similarly, the formula for $\pi_3$ and $\pi_4$ is obtained by applying the trivialization $T\g \rightarrow \g \times T_0\g$ to the element of $T\g$ represented by $\gamma(s,0)^{-1}\gamma(s,t).$
\end{proof}

\begin{corollary}
Let $[c],[d] \in \g,$ and let $\gamma(s,t) = c(s)d(t)c(s)^{-1}d(t)^{-1}.$ Then the trivialization of $T^2G$ maps the element $\xi \in T^2G$ represented by $\gamma$ to
\[
(e,0,0, [s \mapsto [t \mapsto \gamma(s,t)]])
\]
\end{corollary}
\begin{proof}
This follows since $\gamma(s,0) = \gamma(0,t) = e$ for all $s$ and $t.$ 
\end{proof}

\begin{remark}
If $v = [c] \in \g,$ one easily checks that $\gamma_1(s,t) = c(s)$ and $\gamma_2(s,t)=c(t)$ satisfy the conditions of Lemma \ref{lem:curvesinTX}. Let $\iota_1(v)$ and $\iota_2(v)$ denote the corresponding elements of $T^2G.$ This yields canonical maps $\iota_i: \g \rightarrow T^2G,$ which in above trivialization of $T^2G$ correspond to the inclusions of $\g$ as second and third factor of $G \times \g \times \g \times T_0\g,$ respectively. Let $\kappa: G \times G \rightarrow G$ denote the commutator map $\kappa(g,h) = ghg^{-1}h^{-1}.$ Then it follows from Remark \ref{rem:secondDerivative} that $T^2\kappa(\iota_1(v),\iota_2(v))$ is represented by
\[
\gamma(s,t) = \kappa(\gamma_1(s,t),\gamma_2(s,t)) = c(s)d(t)c(s)^{-1}d(t)^{-1}.
\]
\end{remark}

\subsection*{Lie Bracket}

Under the assumption that a natural map $\g \rightarrow T_0\g$ is bijective we can define a bilinear bracket operation $[\cdot, \cdot]$ on $\g.$ Goal of this subsection is to prove that $[\cdot, \cdot]$ is a Lie bracket.

\begin{definition} \label{def:bracket}
a) Suppose that $V$ is a vector space with a Fr\"olicher structure such that scalar multiplication is smooth. We define a map $\Xi:V \rightarrow T_0V$ by sending $v$ to $[s \mapsto sv],$ the line through $v.$ If $V =T_xX \subset TX$ for some Fr\"olicher space $X,$ and $V$ carries the subspace structure, then it is easily checked that $\Xi$ is injective. 

b) Let $G$ be a Fr\"olicher group. The tangent space $\g=T_eG$ inherits a Fr\"olicher structure from $TG.$  Suppose that the map $\Xi:\g \rightarrow T_0\g$ is bijective. Let $v = [c], w= [d] \in \g,$ and let $\xi \in T_0\g$ be represented by $\gamma(s,t) = c(s)d(t)c(s)^{-1}d(t)^{-1}.$ Then we define
\[
[v,w] = \Xi^{-1}(\xi).
\]
\end{definition}

Here is an auxiliary lemma which will be needed later.
\begin{lemma} \label{lem:aux}
Let $X$ and $Y$ be Fr\"olicher spaces, such that the maps $\Xi$ corresponding to $T_xX$ and $T_yY$ are bijective. Let $\phi:X \rightarrow Y$ be smooth. Then
\[
T\phi(\Xi^{-1}(\xi)) = \Xi^{-1}(T^2\phi(\xi)).
\]
\end{lemma}
\begin{proof}
Suppose that $\Xi^{-1}(\xi) \in T_xX$ is represented by $s \mapsto c(s).$ Then $T\phi(\Xi^{-1}(\xi))$ is represented by $\phi \circ c.$ Now we apply $\Xi$ and get an element of $T_0T_xX$ represented by $\gamma(s,t)= \phi(c(st)).$ By Remark \ref{rem:secondDerivative} this map also represents $T^2\phi(\xi),$ which proves our claim.
\end{proof}

Our goal is to prove that $[ \cdot, \cdot]$ is a Lie bracket on $\g.$ To this end, we view elements of $\g$ as derivations of the algebra $F$ of functions on $G.$ It will then suffice to show that $[v,w]$ corresponds to the commutator of the derivations corresponding to $v$ and $w.$

\begin{definition} \label{def:derivation}
We let $\lambda_g: G \rightarrow G$ denote left multiplication by $g \in G.$ If $v = [c] \in \g,$ we define $D_v:F \rightarrow F$ by
\[
D_vf(g) = b(v,[f \circ \lambda_g]) =  (f \circ \lambda_g \circ c)'(0).
\]
\end{definition}

\begin{lemma}
For each $v \in \g,$ the map $D_v$ is a derivation of the $\R$-algebra $F.$ The map $v \mapsto D_v$ is a linear injection of $\g$ into the vector space of derivations of $F.$
\end{lemma}
\begin{proof}
Let $v = [c] \in \g$ and $f,g \in F.$ The vector space $F$ is an algebra under pointwise multiplication of functions. Therefore, $D_v(fg)(x)$ is the derivative of the product $(f \circ \lambda_x \circ c)(g \circ \lambda_x \circ c)$ of functions in $C^\infty(\R,\R),$ and the product rule for derivatives shows that $D_v(fg)(x) = f(x)D_v(g)(x) + g(x)D_v(f)(x).$

Now suppose $v=[c] \neq w=[d]$ in $\g.$ By definition, this means that there is some $f \in F$ with $(f \circ c)'(0) \neq (f \circ d)'(0),$ hence $D_v(f)(e) \neq D_w(f)(e),$ proving injectivity of $v \mapsto D_v.$

Linearity in $v$ follows from $D_v(f)(x) = b(v, [f \circ \lambda_g]).$ 
\end{proof}

Let $[D_v,D_w]$ be the usual commutator of derivations. If 
\begin{equation} \label{eq:comm}
D_{[v,w]} = [D_v,D_w],
\end{equation}
then it follows immediately that $(\g, [\cdot, \cdot])$ is a Lie algebra. The rest of this section is devoted to proving Equation (\ref{eq:comm}).

Let us first compute the left hand side of Equation (\ref{eq:comm}).
\begin{lemma}
Suppose that $\Xi: \g \rightarrow T_0\g$ is an bijection. Let $\xi \in T_0\g$ be represented by $\gamma:\R^2 \rightarrow G.$ If $v = \Xi^{-1}(\xi)$ and $f \in F,$ then
\[
D_v(f)(e) = \frac{\partial^2 (f \circ \gamma)}{\partial s \partial t}(0,0).
\]
\end{lemma}
\begin{proof}
Suppose that $v = [c].$ We have seen that then, $\xi = \Xi(v)$ can be represented by $\rho(s,t) = c(st).$ Let $\bar \gamma(s) = [t \mapsto \gamma(s,t)]$ and $\bar \rho(s) = [t \mapsto c(st)].$ By definition of $T_0\g,$ for every smooth function $\phi:\g \rightarrow \R$ we have 
\begin{equation} \label{eq:derivations}
(\phi \circ \bar \rho)'(0) = (\phi \circ \bar \gamma)'(0).
\end{equation}
We can apply this if $\phi$ is the second component of $Tf:TG \rightarrow \R^2$ for $f \in F,$ to get
\[
\pi_2(Tf(\bar \rho(s))) = (f \circ \rho(s, \cdot))'(0) = \frac{(\partial f \circ \rho)}{\partial t}(s,0),
\]
Hence, with $\phi = \pi_2 \circ Tf,$ the left hand side of Equation (\ref{eq:derivations}) becomes 
\[
\frac{\partial^2 (f \circ \rho)}{\partial s \partial t}(0,0).
\]
We have $\rho(s,t) = c(st),$ so the chain rule yields $\frac{\partial}{\partial t} \rho(s,t) = s(f \circ c)'(st).$ Now we use the product rule to get
\[
\frac{\partial^2 (f \circ \rho)}{\partial s \partial t}(s,t) = (f \circ c)'(st) + s^2(f \circ c)''(st).
\]
We let $(s,t) = (0,0)$ and equate with the right hand side of Equation (\ref{eq:derivations}) to get
\[
(f \circ c)'(0) = \frac{\partial^2 (f \circ \gamma)}{\partial s \partial t}(0,0),
\]
which is what we claimed.
\end{proof}

Now we turn to the right hand side of Equation (\ref{eq:comm}).

\begin{lemma}
If $v=[c], w = [d] \in \g,$ then
\[
D_v(D_wf)(e) = \frac{\partial^2}{\partial s \partial t} \Big|_{(s,t) = (0,0)} f(c(s)d(t)).
\]
\end{lemma}
\begin{proof}
By definition, $D_wf(g) = (f \circ \lambda_g \circ d)'(0).$ Note that $\lambda_e$ is the identity map on $G,$ so that $D_v(D_wf)(e)$ is the derivative of 
\[
s \mapsto D_wf(c(s)) = (f \circ \lambda_{c(s)} \circ d)'(0)
\]
at 0. Now $f(\lambda_{c(s)}(d(t))) = f(c(s)d(t)),$ hence we obtain
\[
D_v(D_wf)(e) = \frac{\partial}{\partial t} \Big|_{t=0} \left(\frac{\partial}{\partial s}\Big|_{s=0} f(c(s)d(t)) \right),
\]
as was claimed.
\end{proof}

Our last step now is to show that left hand side and right hand side agree.

\begin{lemma}
If $v =[c], w = [d] \in \g$ and $\gamma(s,t)= c(s)d(t)c(s)^{-1}d(t)^{-1},$ then
\[
\frac{\partial^2}{\partial s \partial t} \Big|_{(s,t) = (0,0)} f \circ \gamma = \frac{\partial^2}{\partial s \partial t} \Big|_{(s,t) = (0,0)} \left( f(c(s)d(t)) - f(d(s)c(t)) \right)
\]
\end{lemma}
\begin{proof}
First we fix the variable $s$ and use the formulas (\ref{eq:addition}) for inversion and addition in $T_{c(s)}G$ to obtain
\begin{equation*} \begin{split}
\frac{\partial}{\partial t}\Big|_{t=0} (f(c(s)d(t)) - f(d(t)c(s))) &= \frac{\partial}{\partial t}\Big|_{t=0} (f(c(s)d(t)) + f(d(t)^{-1}c(s)) \\
&= \frac{\partial}{\partial t}\Big|_{t=0} (f(c(s)d(t)c(s)^{-1}d(t)^{-1}c(s))) \\
&= \frac{\partial}{\partial t}\Big|_{t=0} (f(\gamma(s,t)c(s)))
\end{split} \end{equation*}
Now fix $t$ and consider the curve $s \mapsto \gamma(s,t)c(s)$ through $e.$ We use the formula for addition in $T_eG$ to get
\[
\frac{\partial}{\partial s}\Big|_{s=0} (f(\gamma(s,t)c(s))) = \frac{\partial}{\partial s}\Big|_{s=0} (f(\gamma(s,t)) +f(c(s))).
\]
Finally we get the desired result by differentiating the last expression with respect to $t.$
\end{proof}

We combine our results to get the following 

\begin{theorem}
Let $G$ be a Fr\"olicher group, $\g = T_eG$ and $\Xi: \g \rightarrow T_0\g$ the map sending $v$ to $[s \mapsto sv].$ If $\Xi$ is bijective, then $[\cdot, \cdot]$ as in Definition \ref{def:bracket} is a Lie bracket. Let $H$ be another Fr\"olicher group for which $\mathfrak h = T_eH$ has the property that $\Xi: \mathfrak h \rightarrow T_0 \mathfrak h$ is bijective. If $\alpha:G \rightarrow H$ is a smooth homomorphism, then  $T\alpha$ is a morphism of Lie algebras.
\end{theorem}
\begin{proof}
We have shown that the map $v \mapsto D_v$ in Definition \ref{def:derivation} is linear, injective and satisfies $D_{[v,w]} = [D_v, D_w].$ This immediately implies that $[\cdot, \cdot]$ is a Lie bracket on $\g.$ Let $\alpha: G \rightarrow H$ be a morphism of Fr\"olicher groups, where $H$ is another group for which $\Xi: \mathfrak h  \rightarrow T_0 \mathfrak h$ is bijective. We claim that $T\alpha$ respects the Lie bracket. Let $v= [c], w =[d] \in \g$ and let $\xi \in T^2G$ be represented by $\gamma(s,t)=c(s)d(t)c(s)^{-1}d(t)^{-1}.$ We have seen in Lemma \ref{lem:aux} that $T\alpha([v,w]) = T\alpha(\Xi^{-1}(\xi)) = \Xi^{-1}(T^2\alpha(\xi)),$ and by Remark \ref{rem:secondDerivative}, $T^2\alpha(\xi)$ is represented by $\alpha \circ \gamma.$ Now $\alpha$ is a homomorphism, therefore $\alpha(\gamma(s,t))$ is the commutator of $\alpha(c(s))$ and $\alpha(d(t)).$ This proves
\[
T\alpha([v,w]) = [T\alpha(v),T\alpha(w)].
\]
\end{proof}

\section{Example} \label{sct:Example}

As an example for the theory developed above, we consider the product group $G=(\R^J, +)$ with its product Fr\"olicher structure as described in  Example \ref{ex:limits}.

The vector spaces $\R^J$ with their product topology play an important role in the theory of pro-Lie groups, see Hofmann and Morris \cite[Appendix 2]{HoMo:07} for a summary of results. Their topological properties depend on the cardinality of $J.$ For instance, $\R^J$ is polish (that is, seperable and completely metrizable) if and only if $J$ is at most countable. It is separable if and only if the cardinality of $J$ is less than or equal to that of $\R.$

Here we will show that under a rather weak assumption on the cardinality of $J,$ the Fr\"olicher spaces $\R^J$ and $T_0\R^J$ are isomorphic. The crucial idea is that functions $f:\R^J \rightarrow \R$ of the Fr\"olicher structure factor through $\R^{J_0}$ with countable $J_0.$ Since $\R^{J_0}$ is a Fr\'echet space, this allows us to regard them as functions of class $C^\infty$ on $\R^{J_0}.$ The assumption made on $J$ is that its cardinality is strictly smaller than that of any \emph{weakly inaccessible cardinal}. We will use the appendix to indicate which cardinals satisfy this assumption. Let us just mention that the existence of weakly inaccessible cardinals can not be proven within the standard axiom system {\bf ZFC}. Therefore, the following lemma is true for most $J$ that naturally occur.

\begin{proposition} \label{prop:Factorization} 
Let $J$ be a set of cardinality smaller than any weakly inaccessible cardinal. Equip $\R^J$ with the product Fr\"olicher structure $(C,F).$ Then $f \in F$ if and only if there is a countable $J_0 \subset J$ such that $f$ factors through a map $f_0: \R^{J_0} \rightarrow \R,$ and the map $f_0$ is of class $C^\infty$ for the product topology on $\R^{J_0}.$
\end{proposition}
\begin{proof}
If $f$ factors  through a $C^\infty$-map $f_0:\R^{J_0} \rightarrow \R,$ it is clearly smooth along smooth curves, hence in $F.$ 
Now conversely, if $f \in F,$ then $f$ is continuous in the $c^\infty$-topology of $\R^J,$ and it was shown in \cite[Example I.4.27]{Kriegl:1997} that $f$ factors as $f_0:\R^{J_0} \rightarrow \R$ with countable $J_0$ (see Appendix for an outline). Now $f_0$ is smooth along smooth curves, and hence by Corollary \ref{cor:Frechet}, it is of class $C^\infty.$
\end{proof}

\begin{theorem} \label{thm:example}
If the cardinality of $J$ is smaller than any weakly inaccessible cardinal, then
\[
T_0 \R^J \cong \R^J.
\]
\end{theorem}
\begin{proof}
Let $C_0$ denote the set of curves through 0. Then we define $\phi: C_0 \rightarrow \R^J$ by $\phi(c)_j = c_j'(0).$ This map is clearly onto. We need to see that it factors through a bijection $T_0\R^J \rightarrow \R^J.$

Claim:  $\ker \phi = \{ c \in C_0 \;|\; \forall \, f \in F: (f \circ c)'(0) = 0 \}.$ 

The inclusion $\supseteq$ is clear, using $f = \pi_j.$ To prove the reverse inclusion, let $f \in F$ and $c \in \ker \phi.$ Then $f$ factors as $f = f_0 \circ \pi,$ where $\pi: \R^J \rightarrow \R^{J_0}$ with countable $J_0 \subset J.$ Hence $(f \circ c)'(0) = (f_0 \circ \pi \circ c)'(0) = d_0f((\pi \circ c)'(0))$ by the chain rule for $C^\infty$-maps. Now $(\pi \circ c)'(0) = 0,$ since otherwise if the $i$th component is non-zero, we use that $\pi_i \circ \pi \circ c = \pi_i \circ c$ to conclude that $(\pi_i \circ c)'(0) \neq 0,$ contradicting the choice of $c.$ Hence $(f \circ c)'(0)= 0.$ 

Claim: The induced bijection $\bar \phi: T_0 \R^J \rightarrow \R^J$ is smooth and has a smooth inverse. 

To see that $\bar \phi$ is smooth, compose it with $\pi_j$ to get $\pi_j \circ \bar \phi ([c]) = (\pi_j \circ c)'(0).$ But this is the second component of $T\pi_j ([c]) \in T\R = \R^2.$ Now $T\pi_j$ is smooth on $T_0\R^J,$ hence $\pi_j \circ \bar \phi$ is smooth. Now consider $\bar \phi^{-1}.$ Given $x=(x_j)_{j \in J} \in \R^J,$ a representative of $\bar \phi^{-1}(x)$ is given by $[d]$ where $d_j(t) = tx_j.$ Let $c \in C$ and $f \in F.$ Then
\[
T_0f \circ \bar \phi^{-1} \circ c: s \mapsto T_0f([t \mapsto (tc_j(s))_{j \in J}]).
\] 
Again we use the fact that $f$ factors through a countable product $\R^{J_0}.$ We identify $[t \mapsto (tc_j(s))_{j \in J_0}] \in T_0\R^{J_0}$ with $(c_j(s))_{j \in J_0} \in \R^{J_0}$ and obtain
\[
(T_0f \circ \bar \phi^{-1} \circ c)(s) = d_0f((c_j(s))_{j \in J_0}),
\]
which is smooth in $s$ as composition of $C^\infty$-functions. 
\end{proof}

\section*{Concluding Remarks}

It would be desirable to verify $T_0\g \cong \g$ for at least the following classes of Fr\"olicher groups:

\begin{itemize}
\item $\mathcal G = C^\infty(M,G),$ in particular for non-compact manifolds $M.$ The difficulty here, as in the case of $\R^J,$ is to determine the class of smooth functions on $\mathcal G.$ By cartesian closedness, $c$ is a curve if and only if $\tilde c:\R \times M \rightarrow G$ is smooth. If $c$ is a curve through $e,$ we get a map $m \mapsto (\tilde c(\cdot,m))'(0) \in \g.$ If one knew the class of functions, one might be able to show that this assignment yields a map $T_e\mathcal G \rightarrow C^\infty(M,\g).$

\item $G = \Diff(M).$ Here, the difficulty is again in finding the class of functions $f: G  \rightarrow \R.$

\item Pro-Lie groups, and as a test case, the $p$-adic solenoid $\mathbb T_p.$ This is the projective limit of the diagram
\[
\mathbb S^1 \overset{\phi}{\leftarrow} \mathbb S^1 \overset{\phi}{\leftarrow} \mathbb S^1 \dots
\]
with $\phi(z)= z^p.$ The projective limit of this diagram in the category of Hausdorff topological groups is a \emph{topological group with Lie algebra} (see \cite{HoMo:07}), and there is a continuous exponential function $exp: \R \rightarrow \mathbb T_p.$ This map is easily seen to be smooth, and this fact can be used to conclude $T_e\mathbb T_p \cong \R$ as a vector space. However, we have not been able to determine the Fr\"olicher structure on this vector space, so that we were not able to show $T_0\R \cong \R.$ A better knowledge of the set $F$ of functions on $\mathbb T_p$ would be necessary.
\end{itemize}

The examples of $\R^J$ and $\mathbb T_p$ seem to indicate that progress is not possible without using a topology on the Fr\"olicher group under consideration. It is not clear to the author whether this has to be done on a case-by-case basis, or if one should use a topology associated to the Fr\"olicher structure. For example, Michor suggested to restrict attention to Fr\"olicher spaces $(X,C,F)$ for which the initial topology with respect to $F$ agrees with the final topology with respect to $C.$

\section*{Appendix}

In this appendix we provide details for the proof of Proposition \ref{prop:Factorization}.

If $J$ is any index set, let $\R^J$ and $2^J$ carry the product topology, and let $\R^J_0$ and $2^J_0$ denote the subsets of elements with countable support. Consider the following two properties, which we think of as properties of the cardinality $\kappa = |J|.$
\begin{itemize}
\item[(A)] If $f:\R^J_0 \rightarrow \R$ is sequentially continuous, then $f$ factors through a subspace $\R^{J_0}$ with countable $J_0 \subset J.$
\item[(B)] If $f:2^J \rightarrow \R$ is sequentially continuous and vanishes on $2^J_0,$ then $f = 0.$
\end{itemize}

Recall that cardinals can be well-ordered and indexed by ordinals: $\aleph_0$ is the first infinite cardinal and for any ordinal $\alpha >0,$ the symbol $\aleph_\alpha$ denotes the $\alpha$th uncountable cardinal.

Mazur proves the following theorem in \cite{Mazur:1952}.
\begin{theorem}
\begin{itemize}
\item $\aleph_0$ satisfies (A) and (B).
\item If $\aleph_\alpha$ satisfies (A) or (B), so does $\aleph_{\alpha+1}.$
\item If the $\kappa_s, s \in S$ and the cardinality $|S|$ satisfy (A) or (B), so does the union $\kappa = \cup_{s \in S} \kappa_s.$
\end{itemize}
\end{theorem}
This says that (A) and (B) hold for all cardinalities smaller than the first weakly inaccessible cardinal.

We now sketch the proof of Kriegl and Michor's Argument. 
They use it to show that the $c^\infty$-topology of $\R^J$ is not completely regular if $|J| > |\R|.$ Note that Hausdorff and completely regular is equivalent to the existence of a uniform structure inducing the topology, hence $\R^J$ does not have a uniform structure inducing the $c^\infty$-topology. In particular, $\R^J$ with the $c^\infty$-topology is not a topological group.

\begin{lemma}
Let the cardinality $|J|$ satisfy (A) and (B) above. Suppose $f:\R^J \rightarrow \R$ is $c^\infty$-continuous. Then $f$ factors through a countable product $\R^{J_0}.$
\end{lemma}
\begin{proof}
If $A \subset J$ and $x \in \R^J,$ let $x_A$ be defined by 
\[
(x_A)_j = \begin{cases} x_j \; \text{if} \;  j \in A \\ 0 \; \text{else}. \end{cases}
\]
Now if $f:\R^J \rightarrow \R$ is continuous, its restriction to $\R^J_0$ is continuous. On $\R^J_0,$ the $c^\infty$-topology agrees with the restriction of the product topology on $\R^J$ (see \cite[Example I.4.12]{Kriegl:1997}). Thus by property (A), the restriction $f|_{\R^J_0}$ depends only on countably many variables $J_0 \subset J,$ that is, $f(x) = f(x_{J_0})$ for all $x \in \R^J_0.$ We claim that $f(x)=f(x_{J_0})$ for all $x \in \R^J.$ To this end, fix $x \in \R^J.$ We view subsets $A \subset J$ as elements of the topological space $2^J,$ and define a function $\phi_x(A) = f(x_A)-f(x_{A \cap J_0})$ from $2^J$ to $\R.$ By the previous argument, $\phi_x(A)=0$ for countable $A.$ We are done if we show that $\phi_x$ is sequentially continuous, since then by property (B) it follows that $\phi_x(J) = 0,$ hence $f(x) = f(x_{J_0}).$ 

In order to prove sequential continuity, let $A_n \rightarrow A$ in $2^J.$ Then one can show that $\{n(x_{A_n} - x_A) | n \in \N \}$ is a bounded subset of the topological vector space $\R^J.$ This implies that $x_{A_n}$ is Mackey-convergent to $x_A,$ hence also convergent in the $c^\infty$-topology (we refer to \cite[Chapter I, Section 4.9]{Kriegl:1997} for details on this point). The map $f$ is $c^\infty$-continuous, and therefore $\phi_x(A_n) \rightarrow \phi_x(A)$ as $n \to \infty.$
\end{proof}


\bibliographystyle{plain}
\bibliography{diss_article}

\end{document}